\documentclass[twoside]{article}
\usepackage{amsfonts}
\usepackage{amsmath,amssymb}
\usepackage{graphicx}
\usepackage[mathlines]{lineno}

\newcommand*\patchAmsMathEnvironmentForLineno[1]{  \expandafter \let \csname old#1\expandafter \endcsname \csname
#1\endcsname
  \expandafter \let \csname oldend#1\expandafter \endcsname \csname
end#1\endcsname
  \renewenvironment{#1}     {\linenomath \csname old#1\endcsname}     {\csname oldend#1\endcsname \endlinenomath}}
\newcommand*\patchBothAmsMathEnvironmentsForLineno[1]{  \patchAmsMathEnvironmentForLineno{#1}  \patchAmsMathEnvironmentForLineno{#1*}}
\patchBothAmsMathEnvironmentsForLineno{equation}
\patchBothAmsMathEnvironmentsForLineno{align}
\patchBothAmsMathEnvironmentsForLineno{flalign}
\patchBothAmsMathEnvironmentsForLineno{alignat}
\patchBothAmsMathEnvironmentsForLineno{gather}
\patchBothAmsMathEnvironmentsForLineno{multline}

\begin{document}

\title{\vspace{-1in}%
\parbox{\linewidth}{\footnotesize\noindent
\textbf{Applied Mathematics E-Notes, 17(2017), ???-???} \copyright \hfill ISSN 1607-2510
\newline
Available free at mirror sites of http://www.math.nthu.edu.tw/$\sim$amen/} 
\vspace{0pt} \\
Generalized Hermite--Hadamard's inequality for functions convex on the coordinates\footnote{This is a preprint of a paper whose final and definite form is published open access in Applied Mathematics E-Notes. See http://www.math.nthu.edu.tw/~amen/ for the final version.}~\thanks{%
Mathematics Subject Classifications:  26A51, 26D15, 26E60, 41A55}}
\date{}
\author{Eze R. Nwaeze\thanks{
Department of Mathematics, Tuskegee University,
Tuskegee, AL 36088, USA.}\\enwaeze@tuskegee.edu\ }
\maketitle

\begin{abstract}
The aim of this paper is to generalize the Hermite--Hadamard inequality for functions convex on the coordinates. Our composite result generalizes the result of Dragomir in \cite{Drag}. Many other interesting inequalities can be derived from our results by choosing different values of $n\in\mathbb{N}.$ Furthermore, we add to the literature a new result for positive functions convex on the coordinates.
\end{abstract}

\section{Introduction}
A function $f:I\rightarrow\mathbb{R},~\emptyset\neq I\subseteq\mathbb{R},$ is said to be convex on the interval $I$ if the inequality $$f(\lambda x + (1-\lambda)y)\leq \lambda f(x)+(1-\lambda)f(y),$$
holds for all $x, y\in I$ and $\lambda\in [0, 1].$
A well celebrated inequality for the class of convex functions is the Hermite--Hadamard's inequality. The inequality states that for any convex function $f:[a, b]\rightarrow\mathbb{R}$ we have

\begin{equation}\label{Had}
f\Big(\frac{a+b}{2}\Big)\leq \frac{1}{b-a}\int_a^b f(x)dx\leq \frac{f(a)+f(b)}{2}.
\end{equation}
The above Hermite--Hadamard inequality on convex functions have been extensively investigated
by a number of authors, see for example the papers \cite{El,Geo,Chen,Drag2} and the references therein.\\

Now for functions convex on the coordinates, we have the following definition.\\

DEFINITION.
A function $f:R=[a, b]\times [c, d]\rightarrow\mathbb{R}$ is said to be convex on the coordinates if the partial mappings
\begin{equation*}
f_{y}:[a,b]\rightarrow \mathbb{R},~f_{y}(u):=f(u,y)~and~f_{x}:[c,d]
\rightarrow \mathbb{R},~f_{x}(v):=f(x,v)
\end{equation*}
defined for all $y\in [c,d]$ and $x\in [a,b],$ are convex. 

Using \eqref{Had}, Dragomir \cite{Drag} proved the following Hadamard's type result for functions, defined on a rectangle, that are convex on the coordinates.\\

THEOREM 1. \label{drag}
Suppose that $f:R=[a, b]\times [c, d]\rightarrow\mathbb{R}$ is convex on the coordinates on $R.$ Then we have the following inequalities:
\begin{align*}
&f\Big(\frac{a+b}{2},\frac{c+d}{2}\Big )\\
&\leq\frac{1}{2}\Bigg[\frac{1}{b-a}\int_a^bf\Big(x,\frac{c+d}{2}\Big) dx+\frac{1}{d-c}\int_c^df\Big(\frac{a+b}{2},y\Big) dy \Bigg]\\
\nonumber
 &\leq \frac{1}{(b-a)(d-c)}\int_a^b\int_c^df(x, y) ~dx dy\\
\nonumber
&\leq \frac{1}{4(b-a)}\int_a^b\Big[f(x,c) +f(x,d)\Big] dx+\frac{1}{4(d-c)}\int_c^d\Big[f(a, y)+f(b, y)\Big] dy\\
&\leq \frac{f(a, c) + f(a, d) +f(b, c) + f(b, d)}{4}.
\end{align*}
The above inequalities are sharp.\\

Many generalizations, extensions, and improvements of the above result are bound in the literature. We invite the interested reader to see the references \cite{bak2,Farid,Nwaeze,Sari,Farid2}. We present here a recent improvement by Bakula \cite{bak}. Specifically, she obtained the following result:\\

THEOREM 2. \label{bakula}
Suppose that $f:R=[a, b]\times [c, d]\rightarrow\mathbb{R}$ is convex on the coordinates on $R.$ Then we have the following inequalities:
\begin{align*}
&f\Big(\frac{a+b}{2},\frac{c+d}{2}\Big )\\
&\leq\frac{1}{2}\Bigg[\frac{1}{b-a}\int_a^bf\Big(x,\frac{c+d}{2}\Big) dx+\frac{1}{d-c}\int_c^df\Big(\frac{a+b}{2},y\Big) dy \Bigg]\\
\nonumber
 &\leq \frac{1}{(b-a)(d-c)}\int_a^b\int_c^df(x, y) ~dx dy\\
\nonumber
&\leq \frac{1}{8(b-a)}\int_a^b\Big[f(x,c) +f(x,d)+2f\Big(x,\frac{c+d}{2}\Big)\Big] dx\\
&~~~~~~~~~~+\frac{1}{8(d-c)}\int_c^d\Big[f(a, y)+f(b, y) +2f\Big(\frac{a+b}{2},y\Big)\Big] dy\\
&\leq \frac{f(a, c) + f(a, d) +f(b, c) + f(b, d)}{16}+\frac{1}{4}f\Big(\frac{a+b}{2},\frac{c+d}{2}\Big)\\
&~~~~~~+ \frac{f\Big(\frac{a+b}{2}, c\Big) +f\Big(\frac{a+b}{2}, d\Big) +f\Big(a,\frac{c+d}{2}\Big)+ f\Big(b,\frac{c+d}{2}\Big)}{8}.
\end{align*}
The above inequalities are sharp.
\smallskip

In this present paper, we present a further generalization of Theorem 1. We also obtain more results in this direction.



\section{Main results}\label{sec3}
In this section, we will present four theorems. The first three results are geared towards the generalization of Theorem 1. The last result is, to the best of our knowledge, completely new. For the proof of our results, we will need the following lemmas in \cite{Alomari1}. The first lemma is  a generalization of the inequalities  in \eqref{Had} above.\\

LEMMA 1. \label{ML1}
Let $F:[a, b]\rightarrow\mathbb{R}$ be a convex function on $[a, b].$ Then the double inequality
\begin{equation}
h\sum_{k=1}^{n} F\left(\frac{x_{k-1}+x_{k}}{2}\right)\leq \int_a^b F(t) dt\leq \frac{h}{2}\left[F(a)+2\sum_{k=1}^{n-1}F(x_k)+F(b)\right]
\end{equation}
holds, where $x_k=a+k\frac{b-a}{n},~~k=0,1,2,\cdots,n;$ with $h=\frac{b-a}{n},~n\in\mathbb{N}.$ The constant 1 in the left-hand side and 1/2 in the right-hand side are the best possible for all $n\in\mathbb{N}.$ If $F$ is concave, then the inequality is reversed.\\

LEMMA 2. \label{ML2}
Let $F:[a, b]\rightarrow\mathbb{R_+}$ be a positive convex function on $[a, b].$ Then the inequality
\begin{equation}
\int_a^b F(x)~dx - (b-a)F(t)\leq \frac{h}{2}\Big[F(a)+F(b)+2\sum_{k=1}^{n-1}F(x_k)\bigg]
\end{equation}
holds for all $t\in[a, b]$, where $x_k=a+k\frac{b-a}{n},~~k=0,1,2,\cdots,n;$ with $h=\frac{b-a}{n},~n\in\mathbb{N}.$ The constant 1/2 in the right-hand side is the best possible, in the sense that it cannot be replaced by a smaller one for all $n\in\mathbb{N}$. If $F$ is concave, then the inequality is reversed.

We now state and prove our first result. \\

THEOREM 3. \label{MR1}
Let $f: R=[a, b]\times [c, d]\rightarrow\mathbb{R}$ be convex on the coordinates on $R$. Then the following inequalities hold

\begin{align*}
&\frac{d-c}{2n}\sum_{k=1}^{n}\int_a^b f\left(x, \frac{y_{k-1}+y_k}{2}\right)dx + \frac{b-a}{2n}\sum_{k=1}^{n}\int_c^d f\left(\frac{x_{k-1}+x_k}{2}, y\right)dy\\
&\leq \int_a^b\int_c^d f(x, y)dxdy\\
&\leq \frac{d-c}{4n}\int_a^b \left[f(x, c) + f(x, d)\right] dx +  \frac{b-a}{4n}\int_c^d \left[f(a, y) + f(b, y)\right] dy\\
&+ \frac{d-c}{2n}\sum_{k=1}^{n-1}\int_a^b f\left(x, y_k\right)dx + \frac{b-a}{2n}\sum_{k=1}^{n-1}\int_c^d f\left(x_k, y\right)dy,
\end{align*}
where $x_k=a+k\frac{b-a}{n}, ~y_k=c+k\frac{d-c}{n}, ~~k=0,1,2,\cdots,n;$ and $n\in\mathbb{N}.$
These inequalities are sharp for each $n.$\\

PROOF.
From the assumption, we have that $g_x:[c, d]\rightarrow\mathbb{R}, ~~g_x(y)=f(x, y)$ is convex on $[c, d]$ for all $x\in [a, b]$. Applying Lemma 1 to the function $g_x$, we obtain

\begin{align*}
\frac{d-c}{n}\sum_{k=1}^n g_x \Big(\frac{y_{k-1}+y_k}{2}\Big)\leq \int_c^d g_x(y)~dy \leq \frac{d-c}{2n}\Big[g_x(c) + 2\sum_{k=1}^{n-1}g_x(y_k)+g_x(d)\Big].
\end{align*}
This implies that 
\begin{align}\label{E1}
\nonumber
\frac{d-c}{n}\sum_{k=1}^n f \Big(x, \frac{y_{k-1}+y_k}{2}\Big)\leq \int_c^d f(x, y)& ~dy \\
&\leq \frac{d-c}{2n}\Big[f(x, c) + f(x, d) + 2\sum_{k=1}^{n-1}f(x, y_k)\Big].
\end{align}
Integrating all sides of \eqref{E1} over the interval $[a, b]$, we get
\begin{align}\label{E2}
\nonumber
\frac{d-c}{n}\sum_{k=1}^n \int_a^b f \Big(x, \frac{y_{k-1}+y_k}{2}\Big)~dx&\leq \int_a^b\int_c^d f(x, y)~dxdy \\
&\leq \frac{d-c}{2n}\Bigg[\int_a^bf(x, c)~dx + \int_a^b f(x, d)~dx + 2\sum_{k=1}^{n-1}\int_a^b f(x, y_k)~dx\Bigg].
\end{align}
Following a similar fashion as outlined above, we have, for the mapping $g_y:[a, b]\rightarrow\mathbb{R}, ~~g_y(x)=f(x, y)$ for all $y\in [c, d]$, the following inequalities
\begin{align}\label{E3}
\nonumber
\frac{b-a}{n}\sum_{k=1}^n f \Big(\frac{x_{k-1}+x_k}{2}, y\Big)\leq \int_a^b f(x, y)& ~dx \\
&\leq \frac{b-a}{2n}\Big[f(a, y) + f(b, y) + 2\sum_{k=1}^{n-1}f(x_k, y)\Big].
\end{align}
Integrating over $[c, d]$, gives
\begin{align}\label{E4}
\nonumber
\frac{b-a}{n}\sum_{k=1}^n \int_c^d f \Big(\frac{x_{k-1}+x_k}{2}, y\Big)~dy&\leq \int_a^b\int_c^d f(x, y)~dxdy \\
&\leq \frac{b-a}{2n}\Bigg[\int_c^d f(a, y)~dy + \int_c^d f(b, y)~dy + 2\sum_{k=1}^{n-1}\int_c^d f(x_k, y)~dy\Bigg].
\end{align}
Adding \eqref{E2} and \eqref{E4}, one gets

\begin{align*}
&\frac{d-c}{2n}\sum_{k=1}^{n}\int_a^b f\left(x, \frac{y_{k-1}+y_k}{2}\right)dx + \frac{b-a}{2n}\sum_{k=1}^{n}\int_c^d f\left(\frac{x_{k-1}+x_k}{2}, y\right)dy\\
&\leq \int_a^b\int_c^d f(x, y)dxdy\\
&\leq \frac{d-c}{4n}\int_a^b \left[f(x, c) + f(x, d)\right] dx +  \frac{b-a}{4n}\int_c^d \left[f(a, y) + f(b, y)\right] dy\\
&+ \frac{d-c}{2n}\sum_{k=1}^{n-1}\int_a^b f\left(x, y_k\right)dx + \frac{b-a}{2n}\sum_{k=1}^{n-1}\int_c^d f\left(x_k, y\right)dy.
\end{align*}
That completes the proof.\\

COROLLARY 1. \label{cor}
Let $f:[0,1]\times[0,1]\rightarrow\mathbb{R}$ be convex on the coordinates. Then the following inequalities hold
\begin{align*}
&\frac{1}{4}\int_0^1 \Big[f\Big(x,\frac{1}{4}\Big) + f\Big(x,\frac{3}{4}\Big)\Big]~dx + \frac{1}{4}\int_0^1 \Big[f\Big(\frac{1}{4},y\Big) + f\Big(\frac{3}{4},y\Big)\Big]~dy\\
&\leq \int_0^1\int_0^1 f(x, y)~dxdy\\
&\leq \frac{1}{8}\int_0^1 \Big[f(x,0) + f(x,1)\Big]~dx + \frac{1}{8}\int_0^1 \Big[f(0,y) + f(1,y)\Big]~dy\\
&+\frac{1}{4}\int_0^1 f\Big(x,\frac{1}{2}\Big)~dx+\frac{1}{4}\int_0^1 f\Big(\frac{1}{2},y\Big)~dy.
\end{align*}

PROOF.
The desired inequalities are obtained by taking $n=2$ in Theorem 3 and observing that $x_0=y_0=0,$~$x_1=y_1=\frac{1}{2},$ and $x_2=y_2=1.$\\

EXAMPLE.
Consider the following function $f:[0,1]\times[0,1]\rightarrow\mathbb{R}$ defined by $f(x,y)=xy.$ The function $f$ is convex on each coordinate. To see this, we have that for $\lambda\in[0,1],$ $f(\lambda x+(1-\lambda)x,y)=f(x,\lambda y+(1-\lambda)y)=\lambda f(x,y)+(1-\lambda)f(x,y).$
It is easy to see that $f$ satisfies the conclusion of  Corollary 1.  In fact, equality is attained.\\

THEOREM 4. \label{MR2}
Under the assumptions of Theorem 3, we have 
\begin{align*}
&\sum_{k=1}^{n}f\left(\frac{a+b}{2}, \frac{y_{k-1}+y_k}{2}\right) + \sum_{k=1}^{n}f\left(\frac{x_{k-1}+x_k}{2}, \frac{c+d}{2}\right) \\
&\leq \frac{n}{d-c}\int_c^d f\left(\frac{a+b}{2}, y\right) dy + \frac{n}{b-a}\int_a^b f\left(x, \frac{c+d}{2}\right) dx. 
\end{align*}
This inequality is sharp for each $n.$

PROOF.
From the first part of \eqref{E1} and \eqref{E3}, we obtain the following inequalities
$$\sum_{k=1}^n f \Big(\frac{a+b}{2}, \frac{y_{k-1}+y_k}{2}\Big)\leq \frac{n}{d-c}\int_c^d f\Big(\frac{a+b}{2}, y\Big) ~dy,$$
and
$$\sum_{k=1}^n f \Big(\frac{x_{k-1}+x_k}{2}, \frac{c+d}{2}\Big)\leq \frac{n}{b-a}\int_a^b f\Big(x, \frac{c+d}{2}\Big)~dx.$$
The desired inequality follows by adding the above two inequalities.\\

If one takes $n=2$ in the above inequality (in Theorem 4), one gets\\

COROLLARY 2.
\begin{align*}
\frac{1}{2}\bigg[f\Big(\frac{1}{2},\frac{1}{4}\Big) &+ f\Big(\frac{1}{2},\frac{3}{4}\Big)+f\Big(\frac{1}{4},\frac{1}{2}\Big)+f\Big(\frac{3}{4},\frac{1}{2}\Big)\bigg]\\
&\leq \int_0^1 f\Big(x, \frac{1}{2}\Big)~dx + \int_0^1 f\Big(\frac{1}{2},y\Big)~dy,
\end{align*}
where $f:[0,1]\times[0,1]\rightarrow\mathbb{R}$ is convex on the coordinates.\\

THEOREM 5. \label{MR3}
Under the assumptions of Theorem 3, we obtain the following inequality
\begin{align*}
&\frac{n}{d-c}\int_c^d \left[f(a, y) + f(b, y)\right]dy +  \frac{n}{b-a}\int_a^b \left[f(x, c) + f(x, d)\right]dx\\
&\leq f(a, c) + f(a, d) + f(b, c) + f(b, d)+ \sum_{k=1}^{n-1}\left[f(a, y_k) + f(b, y_k) + f(x_k, c) + f(x_k, d)\right].
\end{align*}
This inequality is sharp for each $n.$\\

PROOF.
Now using the second part of \eqref{E1} and \eqref{E3}, we get
$$\frac{2n}{d-c}\int_c^d f(a, y)~dy\leq f(a, c) + f(a, d) + 2\sum_{k=1}^{n-1}f(a, y_k),$$
$$\frac{2n}{d-c}\int_c^d f(b, y)~dy\leq f(b, c) + f(b, d) + 2\sum_{k=1}^{n-1}f(b, y_k),$$
$$\frac{2n}{b-a}\int_a^b f(x, c)~dx\leq f(a, c) + f(b, c) + 2\sum_{k=1}^{n-1}f(x_k, c),$$
and
$$\frac{2n}{b-a}\int_a^b f(x, d)~dx\leq f(a, d) + f(b, d) + 2\sum_{k=1}^{n-1}f(x_k, d),$$
which, by addition, gives the desired result.\\

REMARK 1.
By taking $n=1$ in Theorems 3, 4 and 5, and combining the resultant inequalities, we recapture Theorem 1 due to Dragomir.\\

THEOREM 6. \label{MR4}
Let $f: R=[a, b]\times [c, d]\rightarrow\mathbb{R_+}$ be convex on the coordinates on $R$. Then we have the following inequality
\begin{align*}
&\int_a^b\int_c^d f(x, y)~dxdy\\
&\leq \frac{b-a}{4n}\Big[(n+1)\int_c^d f(a,y)~dy+(n+1)\int_c^d f(b,y)~dy+2\sum_{k=1}^{n-1} \int_c^d f(x_k,y)~dy\Big]\\
& +\frac{d-c}{4n}\Big[(n+1)\int_a^b f(x,c)~dx+(n+1)\int_a^b  f(x,d)~dx+2\sum_{k=1}^{n-1} \int_a^b f(x,y_k)~dx\Big],
\end{align*}
where $x_k$ and $y_k$ are defined as in Theorem 3.
This inequality is sharp for each $n.$\\

PROOF.
Applying the inequality of Lemma 2~ to the function $g_y:[a, b]\rightarrow\mathbb{R},~g_y(x)=f(x, y)$ at $x=b,$
\begin{align*}
\int_a^b g_y(x)~dx - (b-a)g_y(b)\leq \frac{b-a}{2n}\Big[g_y(a)+g_y(b)+2\sum_{k=1}^{n-1}g_y(x_k)\Big].
\end{align*}
This implies that 
\begin{align}\label{1}
\int_a^b f(x, y)~dx - (b-a) f(b, y)\leq \frac{b-a}{2n}\Big[f(a,y)+f(b,y)+2\sum_{k=1}^{n-1} f(x_k,y)\Big].
\end{align}
Integrating \eqref{1} over $[c, d],$ we get
\begin{align}\label{2}
\nonumber
&\int_a^b\int_c^d f(x, y)~dxdy\\
&~~~~~~~~\leq \frac{b-a}{2n}\Big[\int_c^d f(a,y)~dy+(1+2n)\int_c^d f(b,y)~dy+2\sum_{k=1}^{n-1} \int_c^d f(x_k,y)~dy\Big].
\end{align}

Applying again Lemma 2 to the mapping $g_y$ at $x=a$ and integrating over $[c, d],$ we have
\begin{align}\label{21}
\nonumber
&\int_a^b\int_c^d f(x, y)~dxdy\\
&~~~~~~~~\leq \frac{b-a}{2n}\Big[(1+2n)\int_c^d f(a,y)~dy+\int_c^d f(b,y)~dy+2\sum_{k=1}^{n-1} \int_c^d f(x_k,y)~dy\Big].
\end{align}
Using \eqref{2} and \eqref{21}, one obtains
\begin{align}\label{22}
\nonumber
&\int_a^b\int_c^d f(x, y)~dxdy\\
&\leq \frac{b-a}{2n}\Big[(n+1)\int_c^d f(a,y)~dy+(n+1)\int_c^d f(b,y)~dy+2\sum_{k=1}^{n-1} \int_c^d f(x_k,y)~dy\Big].
\end{align}

Following a similar fashion for the mapping $g_x:[c, d]\rightarrow\mathbb{R},~g_x(y)=f(x, y)$ at $y=c,$ and $y=d$, and then integrating over $[a, b],$ we get
\begin{align}\label{4}
\nonumber
&\int_a^b\int_c^d f(x, y)~dxdy\\
&\leq \frac{d-c}{2n}\Big[(n+1)\int_a^b f(x,c)~dx+(n+1)\int_a^b  f(x,d)~dx+2\sum_{k=1}^{n-1} \int_a^b f(x,y_k)~dx\Big].
\end{align}
The desired inequality is achieved by adding \eqref{22} and \eqref{4}.\\

REMARK 2.
\begin{enumerate}
\item For $n=1,$ we get
\begin{align*}
&\int_a^b\int_c^d f(x, y)~dxdy\\
&\leq \frac{d-c}{2}\int_a^b \Big[f(x, c) + f(x, d)\Big]~dx + \frac{b-a}{2}\int_c^d \Big[f(a, y) + f(b, y)\Big]~dy.
\end{align*}
\item For $n=2,$ we obtain
\begin{align*}
&\int_a^b\int_c^d f(x, y)~dxdy\\
&\leq \frac{d-c}{8}\int_a^b \Big[3f(x, c) + 2f\Big(x, \frac{c+d}{2}\Big)+ 3f(x, d)\Big]~dx\\
&~~~~~~~~~~~~ + \frac{b-a}{8}\int_c^d \Big[3f(a, y) +2f\Big(\frac{a+b}{2},y\Big)+ 3f(b, y)\Big]~dy.
\end{align*}
\end{enumerate}

\section{Conclusion}
A generalization of the Hermite--Hadamard's inequality for functions convex on the coordinates has been obtained. For different values of $n\in\mathbb{N},$ one obtains loads of new inequalities. In particular, if we take $n=1$ in Theorems 3, 4 and 5, and combine the resultant inequalities, then we recapture Theorem 1 due to Dragomir.  Furthermore, we added to the literature a new result for positive functions convex on the coordinates.\\ \\

\textbf{Acknowledgment.} Many thanks to the anonymous referee for his/her valuable comments.

\end{document}